\documentclass[12pt]{amsart}
\usepackage{latexsym,amssymb,amsfonts,amsmath}
\usepackage{amscd, euscript, MnSymbol}
\usepackage{amscd, euscript}
\usepackage{enumerate}
\usepackage{epsfig,graphicx,graphics} 
\newcommand{\R}{\mathbb{R}}
\newcommand{\K}{\mathbb{K}}
\newcommand{\N}{\mathbb{N}}

\newcommand{\Z}{\mathbb{Z}}
\newcommand{\Q}{\mathbb{Q}}
\newcommand{\gen}[1]{\langle#1\rangle}

\DeclareMathOperator{\Eqv}{Eqv}

\DeclareMathOperator{\Pol}{Pol}
\DeclareMathOperator{\Cong}{Cong}
\newtheorem{definition}{{\bf Definition}}[section]
\newtheorem{definitions}{{\bf Definitions}}[section]
\newtheorem{theorem}[definition]{{\bf Theorem}}

\newtheorem{corollary}[definition]{{\bf Corollary}}
\newtheorem{proposition}[definition]{\noindent {\bf Proposition}}

\newtheorem{claim}[definition]{\noindent {\bf Claim}}

\newtheorem{example}[definition]{\noindent {\bf Example}}
\newtheorem{examples}[definition]{\noindent {\bf Examples}}

\newtheorem{remark}[definition]{\noindent {\bf Remark}}
\newtheorem{problem}[definition]{\noindent {\bf Problem}}
\newtheorem{problems}[definition]{\noindent {\bf Problems}}
\newtheorem{notation}[definition]{{\bf Notation}}

\newtheorem{lemma}[definition]{{\bf Lemma}}
\title{Semirigid systems of three equivalence relations}
\author[C.Delhomm\'e]{Christian Delhomm\'e}
\address{LIM-ERMIT, D\'epartement de Math\'ematiques et Informatique, Universit\'e de La R\'eunion, 15 avenue Ren\'e Cassin - BP 7151 - 97715 Saint-Denis Messag. Cedex 9 FRANCE}
\email[Christian Delhomm\'e]{delhomme@univ-reunion.fr}
\author[M.Miyakawa]{Masahiro Miyakawa}
\address{Tsukuba University of Technology,
4-12-7 Kasuga, Tsukuba,
 Ibaraki 305-8521, Japan} \email {mamiyaka@cs.k.tsukuba-tech.ac.jp}
\author[M.Pouzet] {Maurice Pouzet}
\address{University CDr.de-Bernard  Lyon1,
43, Bd. du 11 Novembre 1918,
69622 Villeurbanne, France et Department of Mathematics and Statistics, The University of Calgary, Calgary, Alberta, Canada}
\email{pouzet@univ-lyon1.fr}
 \author [I.G.Rosenberg]{Ivo G.Rosenberg}
\address{Universit\'{e} de Montr\'{e}al,
C.P.6128,Succ.``Centre-ville'',
Montr\'{e}al P.Q. H3C 3J7,  Canada and Department of Mathematics and Statistics, Masaryk University, Brno, Czeck Republic,}
\email{rosenb@DMS.Umontreal.ca}
\author [H. Tatsumi]{Hisayuki Tatsumi}
\address{Tsukuba University of Technology,
4-12-7 Kasuga, Tsukuba,
 Ibaraki 305-8521, Japan} \email{
tatsumi@cs.k.tsukuba-tech.ac.jp}
\date{\today}

\begin{document}

\begin{abstract}
  A system $\mathcal M$ of equivalence relations on a set $E$  is \emph{semirigid} if  only  the identity  and constant functions preserve all members of $\mathcal M$. We construct  semirigid systems  of three equivalence relations. Our construction leads to the   examples given by Z\'adori in 1983 and to many others and also extends to  some infinite cardinalities. As a consequence,  we show that on every set of at most continuum cardinality distinct from $2$ and $4$ there exists a semirigid system  of three equivalence relations.\end{abstract}
\maketitle
\section{Introduction}
A {\it binary relation} on a set $E$ is a  subset $\rho$ of the cartesian product $E\times E$. We write $x\rho y$ instead of $(x,y)\in \rho$. A map $f:E\rightarrow E$ \emph{preserves} $\rho$ if \begin{equation}
x \rho y \Rightarrow f(x)\rho f(y) \label{eq:1}
\end{equation}
for all $x,y\in E$.
A \emph{binary system} on the set  $E$ is a pair $\mathcal M:=(E, (\rho_i)_{i\in
I})$ where each $\rho_i$ is a binary relation on $E$. An \emph{endomorphism} of $\mathcal M$ is any map $f:E\rightarrow E$ preserving  each $\rho_i$. The identity map on $E$ is an endomorphism of $\mathcal M$. If there is no other endomorphism, $\mathcal M$ is \emph{rigid}.  Provided that the relations $\rho_i$ are reflexive, the constant maps are endomorphisms,  too. We say  $\mathcal M$ is  \emph{semirigid} if the identity map and the constant maps  are the  only endomorphisms.

Rigidity and semirigidity have attracted some attention (eg see \cite {LaP84}, \cite{rosenberg},  \cite{vopenka}). Systems of equivalence relations lead  to prototypes of semirigid systems. Indeed,  as mentioned by R. S. Pierce \cite{Pie68} (Problem 2, p. 38), \emph{if a set $E$ has
at least  three elements, only the constant functions and the identity map preserve all equivalence relations on $E$}. From this,  it follows that:

\begin{lemma}\label{Pierce}
 If a set  $\{\rho_i: i\in I\}$ of equivalence relations generates by means of joins and meets (possibly infinite) the
 lattice $\Eqv(E)$ of equivalences relations on $E$,  then $\mathcal M:=(E, (\rho_i)_{i\in I})$ is semirigid.
\end{lemma}

The converse does not hold.  Indeed, according to Strietz \cite {Stri 77}, if $E$ is finite with at least  four elements,  four equivalences are needed to generate $\Eqv(E)$ by joins and meets whereas, as shown by  Z\'adori  (1983) \cite{Zad83},  for every  set $E$, whose size $\vert E \vert$ is finite and distinct  from $2$ and $4$,  there is a semirigid system made of three equivalence relations. Z\'adori's  result reads as follows.

\begin{theorem}\label{zadori}  Let $A:= \{0, \dots, n-1\}$. The following system of three equivalence
relations $\rho,\sigma,\tau$  on $A$  is  semirigid.

\begin{enumerate}[{}]
\item Case $n=2k+1$, $k\geq 1$:
\begin{eqnarray*}
\rho&=&\{\{0,1,\ldots, k-1\},\{k,\ldots, 2k\}\},\\
\sigma&=&\{\{0, k\},\{1, k+1\}, \ldots,\{k-1, 2k-1\}\},\\
\tau&=&\{\{0, k+1\},\{1, k+2\}, \ldots,\{k-1, 2k\}\}.
\end{eqnarray*}

\item Case $n=2k+2$, $k\geq 2$: 
\begin{eqnarray*}
\rho&=&\{\{0\},\{1,2,\ldots,k\},\{k+1,\ldots, 2k+1\}\},\\
\sigma&=&\{\{0,1,k+1\},\{2,k+2\}, \ldots,\{k,2k\}\},\\
\tau&=&\{\{1,k+2\},\{2,k+3\}, \ldots, \{k-1, 2k\}, \{0,k,2k+1\}\}.
\end{eqnarray*}

\end{enumerate}
\end{theorem}

\begin{figure}[htb]
\begin{center}
\includegraphics[width=5.2in]{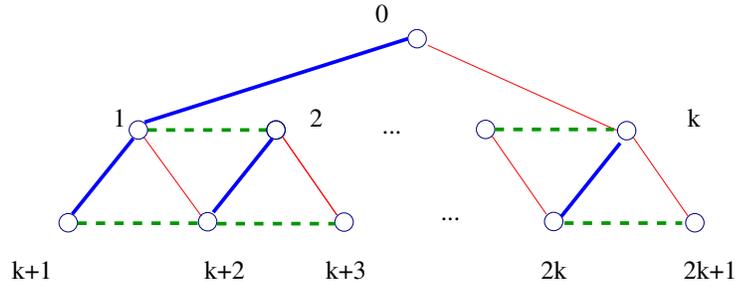}
\end{center}
\caption{One of Z\'adori's examples} \label{Zadori}
\end{figure}

 The case  $n$ even is represented by the colored graph of Figure \ref{Zadori}. In this graph, connected components formed by  single colored edges represent the blocks of a partition into equivalence classes. For $n$ odd, simply delete the node 0 which is located on top of the graph and relabel conveniently the vertices.

In \cite {MPRT} we investigated semirigid systems. In this paper we describe  a general construction of semirigid systems of three equivalence relations  which  includes Z\'adori's. 

\subsection{Results}
Let $\R$ denote  the set of real numbers  and  $\R\times \R$ its cartesian square. We define three equivalence  relations $\simeq_0,\simeq_1$ and $\simeq_2$ on $\R\times \R$. We denote by $p_1$ and $p_2$ the first and second projections from $\R\times\R$ onto $\R$ and  by $p_0: \R\times \R \rightarrow \R$ the map $p_0:=p_1+p_2$. Next, for $i= 0, 1, 2$, we define  $\simeq_i$  as the kernel of $p_{i}$, \emph{i.e.} $u\simeq_iv$ for all $u, v \in \R\times \R$ such that $p_i(u)=p_i(v)$. Finally,  we set $\mathcal R:=(\R\times \R,  (\simeq_0,\simeq_1, \simeq_2))$.

The system $\mathcal R$ is,  by far, not semirigid. But, there are many subsets $C$ of the plane for which the system $\mathcal R\restriction C$ induced by $\mathcal R$ on $C$ is semirigid. In Section \ref{section:abelian} we introduce the notion of \emph{monogenic subset} (see Definitions  \ref {def:monogenic} given in  Subsection \ref{subsection:monogenic}) and with the  simple notion of symmetry in the plane  we prove:  

\begin{theorem} \label{maintheo} If   a finite subset $C$ of $\R\times \R$ is monogenic and has no center of symmetry  in  $\R\times \R$ then $\mathcal R\restriction C$ is semirigid.
\end{theorem}

A simple example, which is at the origin of this result, is the following.
Let $\N$  denotes the set of non-negative integers. For $n\in \N$ set $T_{n}:=\{(i,j) \in \N\times\N: i+j\leq n\}$. Then $T_{n}$ satisfies the hypotheses of Theorem \ref{maintheo}. Hence  $\mathcal R\restriction T_{n}$ is semirigid. This yields an  example of  a semirigid system of three  pairwise isomorphic equivalence relations on a set  having  $\frac{(n+1)(n+2)}{2}$ elements. Now, set $T_{n,2}:=\{(i,j) \in T_{n}: i+j\in \{n-1, n\}\}$ and $T'_{n,2}:=T_{n,2}\cup\{(0,0)\}$. Both sets satisfy the hypotheses of Theorem \ref{maintheo}, hence the induced systems  $\mathcal R\restriction T_{n, 2}$ and $\mathcal R\restriction T'_{n, 2}$ are semirigid. As it is easy to see,  these two examples are isomorphic to those of   Z\'adori.

\begin{figure}[htb]
\begin{center}
\includegraphics[width=5in]{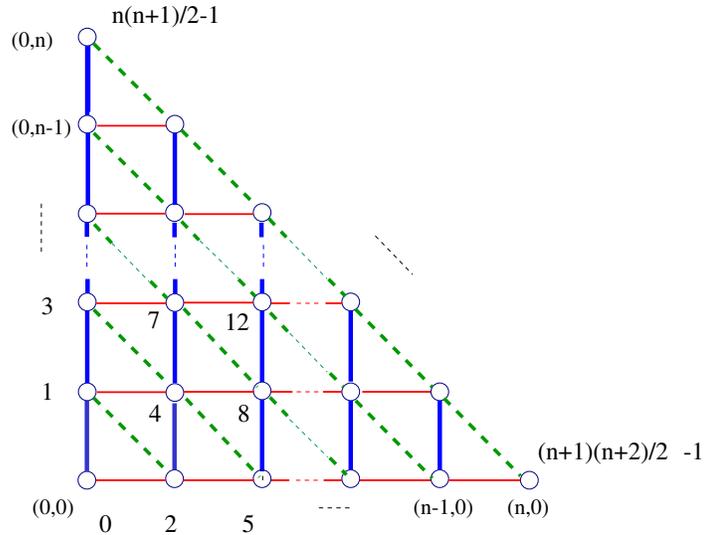}
\end{center}
\caption{The set $T_n$ with the three equivalence relations} \label{3}
\end{figure}

In the vein of Theorem  \ref{maintheo} we prove the following result (see Proposition \ref{prop. semirigid} for a more precise statement). 
\begin{theorem} \label{thm:infinite}For each cardinal $\kappa=3$ or  $4<\kappa \leq 2^{\aleph_0}$, there exists a semirigid system  of three equivalences on a set of cardinality $\kappa$.
\end{theorem}

Our study leads to several questions that we mention  in the  paper. Among these are the following.

\begin{problems} \label{problems}

\begin{enumerate}
\item Does the conclusion of Theorem \ref{thm:infinite} extend to every cardinal $\kappa>2^{\aleph_0}$?

\item A complete description or  a useful characterization of semirigid systems is not known. Find algorithms deciding in reasonable time whether  a system of $\kappa$ equivalence relations on a set of size $n$ is semirigid or not. The algorithmic complexity of the  problem does not seem to be known.
\item Systems of equivalence relations leading to Theorem \ref{maintheo} are subsystems  of systems of equivalence relations associated with three directions in the plane.  Describe  finite semirigid systems  which are embeddable into the system associated with three directions in  the plane. 
\item Describe a test  to decide whether a systems of three  equivalence relations  is  embeddable  into a system   of equivalences associated to three directions in the plane. 
\end{enumerate}
\end{problems}

Concerning item $(3)$ and $(4)$ of Problems \ref{problems} let us mention that the class of systems embeddable into a system of equivalences given by three directions in the plane is closed under
embeddability. This fact leads to he following problems
\begin{problems}
Which are the finite minimal non representable systems?
How many non isomorphic are they? Finitely many?
\end{problems}

\subsection{Links with universal algebra and combinatorics}

To conclude, we indicate  first how our study fits in the context of universal algebra.

  The set of congruences of a universal algebra on  a set $E$ is a basic example  of  a set of equivalence relations.  This  set  ordered by inclusion forms  a lattice, in fact a sublattice of the lattice $\Eqv(E)$.  One of the oldest unsolved problem in  universal algebra is the finite lattice representation problem:
  
\emph{Is  every finite lattice  isomorphic to the congruence lattice of  a finite algebra?}
(see \cite{palfy1, palfy2}). 

An approach to this problem is related to our study. Indeed, for $R$   a subset of $\Eqv(E)$, let $\Pol(R)$, resp. $\Pol^{1} (R)$,  be the set  finitary operations, resp. unary operations,  on $E$ which preserve all members of $R$ and for a set $F$ of finitary operations on $E$, let $\Eqv(F)$ be the set of equivalence relations preserved by  all members of $F$. Ordered by inclusion, the set $\Eqv(F)$ is a lattice, the congruence lattice $\Cong(\mathcal A)$ of the algebra $\mathcal A:=(E, F)$. According to a result of A.Mal'tsev, $\Pol (\Eqv (F))$ is determined by its unary part $\Pol^{1}(\Eqv (F))$, hence the above problem  amounts to prove that every finite lattice $L$ is isomorphic to the congruence lattice of an algebra defined by unary operations.   In general, a sublattice $L$ of $\Eqv(E)$ is not the congruence lattice of an algebra on $E$.  Indeed,  Pierce's result implies that    relational system $\mathcal R:=(E, (\rho_i)_{i\in I})$  made of equivalence relations is semirigid if and only if $\Cong(\mathcal A_{R})=\Eqv(E)$ (where $\mathcal A_{R}:=(E, \Pol^{1}(R))$ and $R:=\{\rho_i: i\in I\}$). 

We recall that $M_3$ is the lattice on five element made of a $3$-element antichain and a top and bottom added. Three equivalences  together with the equality relation  and the full equivalence relation on $E$ form an sublattice of  $\Eqv(E)$  isomorphic to $\mathcal M_3$ if   $(a)$ the meet of any two distinct equivalence relations is the equality relation (this amounts to the fact that the intersection of any pair of equivalence classes belonging to two distinct equivalence relations contains at most one element) and $(b)$ the join of two is the full set $E^2$ of ordered pairs. 

As shown  by Z\'adori \cite{Zad83}:
  
  \begin{theorem}\label{thm:zadori2}
  If $\mathcal R:=(E, (\rho_i)_{i<3})$ made of three equivalences relations is semirigid then  these equivalence relations together with the equality relation  and the full equivalence relation on $E$ form a sublattice $L$ of $\Eqv(E)$, isomorphic to the lattice $\mathcal M_3$.   \end{theorem}

  Provided that $E$ has more than three elements,  the lattice $L$,  mentionned in Theorem \ref{thm:zadori2},   is  distinct from  $\Cong (\mathcal A_R)$ (which is equal to $\Eqv(E)$ since $\mathcal R$ is semirigid) hence is not the congruence lattice of an algebra on $E$.  So our study of semirigidity is somewhat opposite to the representability  question mentioned above. Still, the semirigidity notion appears (briefly) in the thesis of  W. DeMeo \cite{demeo} devoted to the representation of finite lattices.  The sublattices $L$ of $\Eqv(E)$ such that $\Cong(\mathcal A_L)=\Eqv(E)$ (where $\mathcal A_L:=(E,\Pol^{1}(L))$)  are said to be \emph{dense};  the fact that,  as a lattice, $\mathcal M_3$ has a dense representation in every $\Eqv(E)$ with $E$ finite on at least five elements, amounting to Z\'adori's result,  appears  in \cite{demeo} as Proposition 3.3.1 on page 20.  According to  \cite {BKKR69a}, if $R$ is a set of equivalences on a finite set $E$, the congruence lattice $\Cong(\mathcal A_R)$   is made of all  equivalence relations definable by primitive positive formulas  from  $R$ (see \cite{snow} for an easy to read presentation). We  hope that  this result  could help in our problems on  semirigid relations.

Next,  we mention that some systems made of three equivalence relations on the same set  play an important role in the study of latin squares and quasigroups.  Indeed, behind these objects are   $3$-nets, objects usually presented in terms of  incidence structures rather than equivalence relations.  A \emph{$3$-net  of order $n$}  is  a point-line incidence structure made of $n^2$ points and three classes of $n$ lines such that:
\begin{enumerate}[{(i)}]

\item
any two lines from different classes are adjacent;
\item no two lines from the same class are adjacent;
\item any point is incident with exactly one line from each class. 
\end{enumerate}

As the reader will easily observe, this amounts to give three equivalence  relations $\rho_i$, $i<3$, on the set $E$ of points such that each equivalence class of $\rho_i$ intersects in exactly one element each equivalence class of any other $\rho_j$. Let $\rho_i$, $i<3$,  be three equivalence relations on a set $E$ which satisfy this condition (we do not demand on$E$ to be finite); let $E_i$ be the set of equivalence classes of $\rho_i$ and $p_i: E \rightarrow E_i$ be the canonical map. Let $\Phi: E_1\times E_2\rightarrow E_0$ be defined by $\Phi(x,y):=p_0(t)$ where $t$ is the unique element of the intersection $x\cap y$. Then the partial maps $\Phi(x, -): E_2\rightarrow E_0$ and $\Phi  (-,y):E_1\rightarrow E_0$ are bijective for every $x$ and $y$. Conversely,  any map $\Phi$ satisfying these conditions comes from a $3$-net. If the elements of the sets $E_i$ are identified to those of a set $V$, by means of three bijective maps $v_i :V\rightarrow E_i$  we may then define a quasigroup operation $\star$ on $V$ (see subsection \ref{quasigroup}) by setting $a\star b:= v_0^{-1}(\Phi (v_1(a), v_2(b)))$. If $V$ is an $n$-element set, e.g. $V:= \{1, \dots, n\}$, the  $n$ by $n$ matrix whose coefficient $a_{i,j}$ is equal to $i\star j$ is the \emph{table of the quasigroup}. Any table of this form is a \emph{latin square  on the symbols $1, \dots ,n$}. A given $3$-net yields several quasigroups which are said \emph{isotopic}.  

Hundreds of papers and several books, notably \cite{denes-kedwell1, denes-kedwell2},  have been published on latin squares, quasigroups and $3$-nets. We  mention few facts connected to  the problems we are considering. The \emph{graph of 
 a $3$-net} $\mathcal M:= (E, (\rho_i)_{i<3})$ is the  graph $G(\mathcal M)$ with vertex set $E$ and edges the pairs of distinct vertices which belong to some equivalence class of some $\rho_i$. The \emph{graph of a latin square} is the graph of the corresponding $3$-net. Phelps \cite{phelps} has shown that for each integer $n\geq 7$ there is a latin square of order $n$ whose graph has no proper automorphism. No (reflexive) graph with more than a vertex can be semirigid, but it remains to be seen whether the  $3$-nets of Phelps are semirigid. 
 
 \begin{problem}
 
 Does there exist a semirigid $3$-net of order $\kappa 
$ for every cardinal $\kappa$ (finite or not) larger or equal to $7$?

 \end{problem} 

 If $\mathcal M$ is a  $3$-net, the equivalence relations $(\rho_i)_{i<3}$  together with the equality relation  and the full equivalence relation on $E$, form a sublattice  of $\Eqv(E)$, isomorphic to the lattice $M_3$ hence satisfy conditions $(a)$ and $(b)$ beforeTheorem \ref{thm:zadori2} above.  If $\mathcal M$ is induced by a $3$-net on some subset then it satisfies condition $(a)$. Conversely  a system $\mathcal M$ satisfying condition $(a)$   can be extended to a system $\mathcal M'$ which forms a $3$-net. Furthermore, if $E$ is finite one can find $\mathcal M'$ with order at most $2\times \vert E\vert$ (Proposition \ref{prop:extending}).  In the finite case, this  is a rewording of a result about extension of partial latin squares due to Evans \cite{evans}. In  the infinite case, this follows from Compactness theorem of first order logic. Despite some possible confusion, we   call \emph{partial $3$-net} any  system induced by some $3$-net. 

The addition  on the set $\R$ of real numbers is a prototypal example of  a quasigroup operation. The $3$-net associated with this operation is made of the three equivalence relations  $\simeq_0,\simeq_1$ and $\simeq_2$ on $\R\times \R$.  Our  results  are about  partial $3$-nets induced by this  $3$-net.

This paper is divided into two more sections. In the first section, we recall the definitions that we need. We indicate that systems of equivalence relations are  pre-ultrametric spaces, a fact which may also motivate the study of the semirigidity of these systems. We give a representation theorem  for reduced systems of equivalence relations (Theorem \ref{theo1}). We end the section with some examples of semirigid systems; Theorem \ref {theo3} contains, with an other proof,  the fact that $\mathcal R\restriction T_{n}$ is semirigid. Results announced in the introduction are proved in  the last section that is self-contained up to subsection \ref{subsection:self-contained}
.

 \section{Systems of equivalence relations and their representations}
\subsection{Basic definitions}\label{subsection:self-contained}
An {\it equivalence relation} on a set $E$  is a binary reflexive,
symmetric and transitive relation on $E$. A {\it system of equivalence relations} on  $E$ is a pair $\mathcal M:=(E, (\rho_i)_{i\in
I})$ where each $\rho_i$ is an equivalence relation on $E$. If $F$ is a
subset of $E$, the  \emph{restriction} of $\mathcal M$ to $F$ is $\mathcal M\restriction F:=(F, ((F\times F) \cap
\rho_i)_{i\in I})$. If $\mathcal M':=(E', (\rho'_i)_{i\in I})$ is another system of equivalence relations on $E'$, a
\emph{homomorphism} from $\mathcal M$ to $\mathcal M'$ is a map $f: E \rightarrow E'$ such that  for all $x,y\in E, i\in I$:

\begin{equation}
x \rho_iy\;  \text {implies}\;  f(x)\rho'_if(y).\label{eq:1}
\end{equation}

If the map $f$ is injective and he implication in    (\ref{eq:1}) is an equivalence then $f$ is an
\emph{embedding};  if furthermore $f$ is bijective, $f$ is an \emph{isomorphism} and  the systems $\mathcal M$ and  $\mathcal M'$ are \emph{isomorphic}.
A homomorphism of $\mathcal M$ into $\mathcal M$ is called an \emph{endomorphism} of $\mathcal M$; the set of endomorphisms of $\mathcal M$ is denoted by $End(\mathcal M)$.

A system of equivalence relations  $\mathcal M:= (E, (\rho_i)_{i\in I})$ is {\it reduced} if $\cap_{i\in I}\rho_i=
\Delta_E$, where $\Delta_E :=\{(x,x):x\in E\}$.  Note that \emph{if $\mathcal M$ is semirigid then it is reduced}. Indeed, supposing  that $Z:=\bigcap_{i\in I}\rho_i$ is distinct from $\Delta_E$,  one choose some $(x,y)\in Z\setminus \Delta_E $. Then the map $f$ defined by $f(x)=y$ and $f(z)=x$ for $z\not = x$ is an endomorphism of $\mathcal M$. Hence $\mathcal M$ is not semirigid.

\subsection{Systems of equivalence relations and ultrametric spaces}
Propositions \ref{prop:ultra1} and  \ref{prop:ultra2} of this subsection express the fact that systems of equivalence relations  and pre-ultrametric spaces are two faces of the same coin. The study of metric spaces with distance values in a Boolean algebra  first appears in Blumenthal \cite {blum}.  A general study of distances with values in an ordered set  is in \cite{pouz-rose}. Ultrametric spaces with   distance values in an ordered set have been studied by Priess-Crampe and Ribenboim in  several papers, e.g.  \cite{ribenboim1}, \cite{ribenboim2},\cite{ribenboim3}.

A \emph{join-semilattice} is an ordered set in which two arbitrary elements $x$ and $y$  have  a join, denoted by $x\vee y$,  defined as the least element of the set of common upper bounds of $x$ and $y$.

Let $V$ be  a join-semilattice with a least element, denoted by $0$.
A \emph{pre-ultrametric space} over $V$ is a pair $\mathcal D:=(E,d)$ where $d$ is a map from $E\times E$ into $V$ such that for all $x,y,z \in E$:
\begin{equation} \label{eq:ultra1} 
d(x,x)=0,\; d(x,y)=d(y,x)  \text{~and } d(x,y)\leq d(x,z)\vee d(z, y).
\end{equation}

\noindent The map $d$ is an \emph{ultrametric distance} over $V$ and $\mathcal D$ is an \emph{ultrametric space} over $V$ if $\mathcal D$ is a pre-ultrametric space and $d$ satisfies \emph{the separation axiom}:
\begin{equation} \label{eq:ultra2}
d(x,y)=0\; \text{implies} \;  x=y .
\end{equation}

Given a set $I$ let $\powerset (I)$ be the power set of $I$. Then $\powerset (I)$, ordered by inclusion, is a join-semilattice (in fact a complete Boolean algebra) in which the join is the union, and  $0$  the empty set.
 \begin {proposition}\label{prop:ultra1}
Let $\mathcal M:=(E, (\rho_i)_{i\in
I})$ be a system of equivalence relations. For $x,y\in E$, set $d_{\mathcal M}(x,y):=\{i\in I: (x,y)\not \in \rho_i\}$. Then   the pair $U_{\mathcal M}:=(E, d_{\mathcal M})$ is a pre-ultrametric space over $\powerset (I)$.

\noindent Conversely, let $\mathcal D:=(E,d)$ a pre-ultrametric space over $\powerset (I)$. For every $i\in I$ set $\rho_i:=\{(x,y)\in E\times E: i\not \in d(x,y) \}$ and let  $\mathcal M:=(E, (\rho_i)_{i\in
I})$. Then $\rho_i$ is an equivalence relation on $E$ and $d_{\mathcal M}=d$.

\noindent Furthermore,  $U_{\mathcal M}$ is an ultrametric space if and only if $\mathcal M$ is reduced. \end{proposition}

For a join-semilattice $V$  with a $0$ and for two pre-ultrametric spaces $\mathcal D:=(E,d)$ and $\mathcal D':=(E',d')$  over $V$,  a \emph{non-expansive mapping} from $\mathcal D$ to $\mathcal D'$ is any map $f:E\rightarrow  E'$ such that for all $ x,y\in E$:
\begin{equation}
d'(f(x),f(y))\leq d(x,y). \end{equation}
\begin{proposition} \label{prop:ultra2} Let $\mathcal M:=(E, (\rho_i)_{i\in
I})$ and $\mathcal M':=(E', (\rho'_i)_{i\in
I})$ be two systems of equivalence relations. A map $f:E\rightarrow E'$ is a  homomorphism from $\mathcal M$ into $\mathcal M'$ if and only if $f$ is a non-expansive mapping from $U_{\mathcal M}$ into $U_{\mathcal M'}$.
\end{proposition}
Propositions \ref{prop:ultra1} and \ref{prop:ultra2} are immediate and the proofs are left to the reader.
Still, this  suggests to study ultrametric  spaces (in the ordinary sense) which are rigid with respect to  non-expansive mappings. But notice  that a non-trivial ultrametric space could not  be semirigid. Indeed:
\begin{lemma}  An ultrametric space $\mathcal D:=(E, d)$ over a  linearly ordered set  $V$ having at least two elements is not semirigid.
\end{lemma}
\begin{proof} If   all non-zero distances between elements of $\mathcal D$  are equal then every permutation of $V$ is an endomorphism of $D$, hence $\mathcal D$ is  non semirigid. Thus assume that there are  $r, r'$ in $V$  and $x,y, x',y'\in E$ such that  $0<r=d(x,y)<r'=d(x',y')$. Let  $B(x, r):=\{z\in E: d(x,z)\leq r\}$. Then $B(x, r)\not =E$. Indeed, if  $z,z'\in B(x, r)$ then $d(z,z')\leq d(z,x)\vee  d(x, z')\leq r\vee r=r$. Hence, $\{x',y'\}\not \subseteq B(x, r)$.  Let $f$ be the map defined by setting $f(z):=x$ if $z\in B(x, r)$ and $f(z)=z$ otherwise. Then $f$ is non-expansive.  Indeed, let $z,z'\in E$. Suppose $z\in  B(x, r)$ and  $z' \notin B(x, r)$.  Since $V$ is totally ordered, $d(x,z')>r$. Since $d(x, z')\leq d(x,z)\vee d(z,z')\leq r\vee d(z,z')$, we have $d(z, z')\not \leq r$. Since $V$ is totally ordered, we also have $d(x,z')\leq d(z,z')$. This  amounts to $d(f(z),f(z'))\leq d(z,z')$. The other cases lead trivially to the same  inequality.  Hence,  $f$ is non-expansive.
\end{proof}

\begin{problem}Beyond the fact that it is not linearly ordered, what can be said about the order structure of the set of values of the distances of a non-trivial semirigid ultrametric space?
\end{problem}

Contrary to the case of ultrametric spaces, there exist  metric spaces (in the ordinary sense) which are rigid with respect to  non-expansive maps. In fact there are metrizable topological  spaces which are rigid with respect to  continuous maps. Indeed, as shown by de Groot \cite{degroot}, there exist connected and locally connected subsets of the plane which are semirigid with respect to  continuous maps (see Corollary 3 in \cite {degroot}).

\begin{problem} Describe the semirigid metric spaces with respect to  non-expansive mappings.
\end{problem}

\subsection{Representation of  systems of equivalence relations} We may also represent systems of equivalence relations on finite sets by matrices.
Let $A:= (a_{ij})_{i=1, \ldots ,m \atop {j=1, \ldots ,n}}$ be an
$m$ by $n$ matrix with non-negative integer coefficients.
For each $j\in \{1, \dots n\}$ let $\rho_j$ be the equivalence relation defined on $\{1, \dots, m\}$ by
$ i\rho_ji' \; \text{ if }\;  a_{ij} = a_{i'j}.$ Then $\mathcal M(A):=(\{1, \dots, m\}, (\rho_j)_{j=1, \ldots ,n})$ is a system of equivalence relations. Conversely:
\begin{lemma}\label{lem} Every system of $n$ equivalence relations on an $\mathcal M$-element set  is isomorphic to a system  $\mathcal M(A)$ for some $m$ by $n$ matrix $A$ with non-negative integer coefficients.
\end{lemma}
\begin{proof}
 Let $\mathcal M$ be a system of $n$ equivalence relations on an $m$-element set. Without loss of generality we may suppose that the $m$-element set is  $\{1, \dots, m\}$ and hence that $\mathcal M=(\{1, \dots, m\}, (\rho_j)_{j=1, \ldots ,n})$. For $i\in \{1, \dots, m\}$ and $j\in \{1, \ldots ,n\}$ set
$\rho_j(i) := \{i'\in \{1, \dots, m\}: i\rho_j i'\}$.
Let $\varphi$ be a one-to-one map from the set $\{ \rho_j(i): i\in \{1, \dots, m\}, j=1, \ldots ,n\}$ into the set of non-negative integers. Set $a_{ij}:=\varphi(\rho_j(i))$.
\end{proof}
\begin{problem} What is the complexity of the problem: decide whether or not  the system of equivalence relations associated with an $m$ by $n$ matrix is semirigid?
\end{problem}

For reduced systems,  Lemma \ref{lem} has a somewhat simpler form.

Let $\mathcal M:=(E, (\rho_i)_{i\in I})$ be a sytem of equivalence relations. For $x\in E$ and for $i\in I$ let 
$\rho_i(x) := \{y\in E: x \rho_i y\}$ be the equivalence class of $x$ and let 
$E_i:= \{\rho_i(x): x\in E\}$ be the set of equivalence classes.
Let $\overline E:=(E_i)_{i\in I}$ be a family of sets and let $\prod_{i\in I}E_i$ be their cartesian product.  For  $x\in \prod_{i\in I}E_i$ denote by $x_i$ its $i$-th projection, hence $x=(x_i)_{i\in I}$. Let ${\mathcal M}(\overline E):=(\prod_{i\in I}E_i, (\sim_i)_{i\in I})$ where $x \sim_i y$ if $x_i=y_i$. Clearly this system is reduced. If all $E_i$ are equal to the  same set $W$, we denote the system by ${\mathcal M}(W,I)$.

\begin{theorem}\label{theo1} Every reduced system can be embedded into a system of the form ${\mathcal M}(V)$.
\end{theorem}

\proof Let $\mathcal M:=(E, (\rho_i)_{i\in I})$. Let $\overline E:=(E_i)_{i\in I}$ where $E_i:= \{\rho_i(x): x\in E\}$ is the set of equivalence classes of $\rho_i$ and let
 $\rho: E\rightarrow \prod_{i\in I} E_i$ be defined by setting $\rho(x):=(\rho_i(x))_{i\in
I}$. Then, the   map $\rho$  is an embedding from $\mathcal M$ into ${\mathcal M}(\overline E)$:
First, $\rho$ is one-to-one. Indeed,
suppose that $\rho(x)=\rho(y)$ for some $x,y\in E$. Then $\rho(x)_i=\rho(y)_i$, that is $\rho_i(x)=\rho_i(y)$,   for all $i\in I$. Since $\mathcal M$ is reduced,
$x=y$. Next, $f$ is  an homomorphism of $\mathcal M$ into ${\mathcal M}(\overline E)$.
Suppose that $x\rho_i y$ for some  $x,y$ and $i\in I$.   Then $\rho_i(x)=\rho_i(y)$, that is $\rho(x)_i=\rho(y)_i$, which amounts  to the fact that $\rho(x)
\sim_i \rho(y)$. The proof of the converse is similar. \endproof

\begin{corollary} \label {cor1} Every reduced system can be embedded into a system of
 the form ${\mathcal M}(W,I)$, where ${\mathcal M}(W,I)=(W^I, (\sim_i)_{i\in I})$.
\end{corollary}

\proof ${\mathcal M}(\overline E)$ is embeddable into ${\mathcal M}(W,I)$, where
$W:=\bigsqcup_{i\in I}E_i$.
Indeed,  to $x\in \overline E:= \prod_{i\in I} E_i$ associate $\varphi(x)\in W^I$ defined by
$\varphi(x)_i=x_i$. \endproof

Let $\mathcal M:=(E, (\rho_i)_{i\in
I})$ be a system of equivalence relations and  $f$ be  a selfmap of $E$. It is immediate to see that   $f$ is an endomorphism of $\mathcal M$ if and only if for each $i\in I$, the map $f$  induces a selfmap  $f_i$ on the set $E_i$ of equivalence classes of $\rho_i$.  Furthermore, if $\mathcal M$ is reduced, the family $(f_i)_{i\in I}$ determines $f$. Theorem \ref{theo2} below expresses  this fact and characterizes  families  $(f_i)_{i\in I}\in E_i^{E_i}$ that come  from an endomorphism.

\begin{notation}Let $ E\subseteq W^I$ and for $i\in I$ denote  $E_i$  the image  of the $i$-th projection of $E$ in  $W$, that is,  with our notations,  $E_i:=\{x_i: x\in E\}$. Note that $E\subseteq \Pi_{i\in I} E_i$. For  $a\in
\Pi_{i\in I}W^{E_i}$, denote  by $\hat a$ the map from $\Pi_{i\in I} E_i$ into $W^I$ defined by setting:
\begin{equation}\hat a((x_i)_{i\in I}):=(a_i(x_i))_{i\in I}  \;\text{for all}\; x\in \Pi_{i\in I}E_i.
\end{equation}
\end{notation}
\begin{theorem}\label{theo2} Let $ E\subseteq W^I$.  A map   $f:
 E\rightarrow W^I$ is a homomorphism  of ${\mathcal M}(W,I)\restriction E$ into ${\mathcal M}(W,I)$
if and only if $f=\hat a\restriction E$ for some $a\in
\Pi_{i\in I}W^{E_i}$. Moreover, if  such  an element $a$ exists it is unique.
\end{theorem}

\proof
Suppose that $f$ is a  homomorphism. Let $i\in I$ and let   $t\in E_i$. Let $x,x'\in E$ such that $x_i=x_i'=t$. Then
$x \sim_i x'$. Since $f$ is a  homomorphism,  $f(x)\sim_i f(x')$, amounting to
$f(x)_i=f(x')_i$. Set
$a_i(t):=f(x)_i$. Clearly $f(x)=(a_i(x_i))_{i\in
I}$. Set  $a:=(a_i)_{i\in I}$.  We have $f=\hat a\restriction E$. The converse is immediate.  \endproof

\begin{corollary}\label{cor3}Let $ E\subseteq W^I$.  Then    ${\mathcal M}(W,I)\restriction E$ is semirigid if and only if for every $a\in \prod_{i\in I}E_{i}^{E_{i}}$
 such that for every $x\in E$,
 \begin{equation}\label{eq:3}(a_i(x_i))_{i\in I}\in E
 \end{equation}
then  either each $a_i$ is a constant function or each $a_i$ is the identity function on $E_i$.
\end{corollary}
\begin{corollary} There is a bijective correspondence between $End(\mathcal M(W,I))$ and $(W^W)^I$.
\end{corollary}

\begin{proof} Let $a\in (W^W)^I$.  Then $\hat a$  is an endomorphism of $\mathcal M(W,I)$.  Theorem \ref{theo2} ensures that every endomorphism has this form. \end{proof}

\begin{example} Let $W:=\{0,1\}$. Then $\mathcal M(W,I)$ has  $4^{\vert I\vert}$ endomorphisms. Each  endomorphism  $f$ of $\mathcal M(W,I)$ can be expressed
as:

$$ f(x)=(\alpha_ix_i+\beta_i)_{i\in I}$$

\noindent where $\alpha_i\in\{0,1\}$, $\beta_i\in\{0,1\}$, and the sum is modulo 2.
\end{example}

\subsection{Examples of semirigid systems}
\begin{theorem}\label{theo3}Let $W:= \N$, $n\in \N$ and let $I$ be  a  set (not necessarily finite) with at least three elements. Set:
$$E:=\{(x_i)_{i\in I}\in \N^{I} : \sum_{i\in I}
x_i=n\}.$$
Then ${\mathcal M}(\N,I)\restriction E$ is semirigid.
\end{theorem}
\proof Observe first that  $E_i=\{0,\ldots,n\}$ for each $i\in I$.
Apply Corollary \ref{cor3}. Let $a\in \prod_{i\in I}E_{i}^{E_{i}}$
be  such that  Condition  (\ref{eq:3}) of Corollary \ref{cor3} holds.

\begin{claim} \label{claim1}If $i,j\in I$ with $i\not =j$ then for all $t<n$: 

$$a_i(t+1)-a_i(t) = a_j(1)-a_j(0).$$
\end{claim}
\noindent Indeed, let $k$ be distinct from $i$ and $j$.  Let $x:= (x_l)_{l\in I}$ and  $x':=(x'_l)_{l\in I}$  be defined by setting $x_l=x'_l=0$ if $l\in I \setminus \{i,j,k\}$,  $x_i=x'_i+1=t+1$, $x_j=x'_j-1=0$ and $x_k=x'_k=n-t-1$. Then  $x, x'\in E$. Since Condition (\ref{eq:3}) holds,  $(a_l(x_l))_{l\in I}$ and $(a_l(x'_l))_{l\in I}$ belong to $E$,  that is  $\sum_{l\in I}a_l(x_l)=\sum_{l\in I}a_l(x'_l)=n$. Thus $\sum_{l\in I}(a_l(x_l)-a_l(x'_l))= a_i(t+1)-a_i(t)+ a_j(0)-a_j(1)=0$,  proving our claim.

\begin{claim} \label{claim2} $a_i(t)=a_i(0)+ t(a_j(1)-a_j(0))$ for all $t\leq n$.
\end{claim}
\noindent Indeed, this equality holds for $t=0$; for larger $t$ apply induction and Claim \ref{claim1}.
\begin{claim}\label{claim3}
The value $\delta:=a_i(1)-a_i(0)$ is independent of $i$.
\end{claim}
\noindent Indeed, from Claim \ref{claim1}  we obtain $a_i(1)-a_i(0)= a_j(1)-a_j(0)$.

If $\delta =0$ then according to Claim \ref{claim2}, we  have that  $a_i(t)=a_i(0)$ for all $i$ and hence all maps  $a_i$ are
constant, proving that $a$ is constant. Thus, we may suppose   $\delta \not= 0$.
\begin{claim}\label{claim4} $\delta =1$ and  $a_i(0)=0$ for all $i\in I$.
\end{claim}
Indeed,   from Claim \ref{claim2}, we have that   $a_i(n)-a_i(0)=n\delta$.  Since  $a_i(n), a_i(0)\in E_i=\{0,\dots, n\}$,  it follows that  $-1\leq \delta  \leq 1$.  Suppose  that $\delta=-1$.  Then $a_i(0)=n$ and
$a_i(n)=0$ for each $i\in I$. Since $\vert I\vert \geq 3$,  the image $a(x)$ of the sequence $x$ such that $x_i=n$ and $x_l=0$ for $l\not = i$  satisfies $a(x)_i:=a_i(x_i)=a_i(n)=0$ and $a(x)_l:=a_i(x_l)=a_l(0)=n$. Since $\vert I\vert \geq 3$, this image is not  in $E$. This case is then impossible. The
only remaining possibility is $\delta=1$. In this case, $a_i(t)=t$ for all $i$ and
$t$. According to Corollary \ref{cor3},  ${\mathcal M}(\N,I)\restriction E$ is semirigid. \endproof

\subsubsection {Two special cases}

1) Set $I:=\{0,1,2\}$. Then, the map $f$ from $E$ into $T_n:=\{(i,j) \in \N\times\N: i+j\leq n\}$ defined by setting $f(x_0,x_1,x_2):=(x_1,x_2)$ for all $(x_0,x_1,x_2)\in E$ a  is an isomorphism from ${\mathcal M}(\N,I)\restriction E$ onto $\mathcal R\restriction T_{n}$. Hence, $\mathcal R\restriction T_{n}$ is semirigid.

2) Let $n=1$ and $\vert I\vert \geq 3$. In this case, $E$ is the set of characteristic functions of singletons, that is functions $\chi_{\{i\}}$ taking value $1$ on $i$ and $0$ on elements $i'$ distinct from $i$. For each $i$, clearly $\sim_i$ is the  equivalence relation which
puts $\chi_{\{i\}}$ in one class and all other characteristic functions in an other class. Hence, this system of equivalence relations is isomorphic to $\mathcal M:= (I, (\rho_i)_{i\in I})$, where $j\rho_i k$ if either $i=j=k$ or $j,k\in I\setminus \{i\}$.  The fact that the system $\mathcal M$ is
semirigid follows directly  from Lemma \ref{Pierce}. Indeed, for all  $i\not= j$, clearly $\bigcap\{\rho_k: k \in I\setminus \{i, j\}\}$
is a co-atom in the lattice of equivalence relations, and all co-atoms are
obtained this way. Thus the $\rho_i$'s generate the lattice of equivalence relations on $I$.

\section{Systems of three equivalences embedded into 	algebraic structures}\label{section:abelian}

The simplest  systems of three equivalence relations are those associated to groups. Their endomorphisms are described in  Theorem \ref{cor: additivemaps}. These systems are not semirigid, but some induced systems are. In order to  obtain some examples of semirigid systems we describe in Lemma \ref{lem:additive} the endomorphisms defined on induced systems.

 Let $G$ be a  group, the operation being denoted by $\cdot$ and the neutral element $1$. Recall that a \emph{group homomorphism} is any map $h$ from $G$ into $G$ such that  $h(x\cdot y)=h(x)\cdot h(y)$ for all $x,y\in A$. Set $E:=G\times G$. Denote by $p_1$ and $p_2$ the first and second projections from $E$ onto $G$ and  let $p_0:E\rightarrow G$ be defined by $p_0(x,y):=x\cdot y$ (hence $p_0=p_1\cdot p_2$). For $i<3$, let  $\simeq_i$ denote the equivalence relation on $E$ defined for all $u, v \in E$  by setting $u\simeq_iv$ if $p_i(u)=p_i(v)$.  Finally, set $\mathcal M:=(E, (\simeq_i)_{i<3})$ and let $C$ be a subset of $E$.

 \begin{lemma}\label{lem:additive} A map   $f: C\rightarrow E$ is  a homomorphism from $\mathcal M\restriction C$ into $\mathcal M$ if and only if for each $i<3$ there is a map $h_i:p_i(C)\rightarrow G$ such that $p_i(f(x,y))=h_i(p_i(x,y))$ for all $(x,y)\in C$.
In particular, these three maps satisfy:
\begin{equation}\label{eq:additivity}
h_1(x)\cdot h_2(y)=h_0(x\cdot y)
\end{equation}
for all $(x,y)\in C$.

 Furthermore, if  $f$ fixes $(1,1)$, then  $h_1, h_2$ and $h_0$ coincide on $\check C:=\{x\in A: (1,x)\;\text{and}\;  (x,1)\in C\}$  and their common value $h$ satisfies:  \begin{equation}\label{eq:homo}
h (x\cdot y)=h(x)\cdot h(y)
\end{equation}
for all  $x,y$ such that $(x,y)\in C\cap (\check C\times \check C)$ and $x\cdot y\in \check C$. \end{lemma}

\begin{proof}The first part of the lemma holds if $G$ is a \emph{magma} (rather than a group), that is just a set equipped with a binary operation; the second part holds if in addition  this operation  has a neutral element, that is  an element, denoted by $1$,  and necessarily unique,  such that $e\cdot x=x\cdot e=x$ for all $x\in G$. Since each $p_i$ identifies with the canonical projection onto the quotient of $\mathcal M/\simeq_i$, a homomorphism $f$ from $\mathcal M\restriction C$ into $\mathcal M$ induces three maps $h_i: p_i(C) \rightarrow  G$, $i < 3$ (cf. the proof of Theorem \ref{theo2}). The fact that  Equation (\ref{eq:additivity}) holds follows readily from the definitions. Indeed, 
 for $(x,y)\in C$, we have:    $$h_0(x\cdot y)= h_0(p_0(x,y))=p_1(f(x,y))\cdot p_2 (f(x,y)), $$
$$p_1(f(x,y))= h_1(p_1(x,y))=h_1(x)\;  \text{and}  \; p_2(f(x,y))= h_2(p_2(x,y))=h_2(x).
$$  

For the second part, let $x\in \check C$. The relation  $(1,1)\simeq (x,1)$ yields $(1,1)=f(1,1) \simeq f(x,1)$ hence $f(x, 1)=  (h_1(x), 1)$; similarly,  $f(1,x)= (1, h_2(x))$. The relation  $(x,1)\simeq_0(1, x)$ yields    $f(x,1)\simeq_0 f(1,x)$ thus $h_1(x)=h_2(x)$. From $f(1,1)=(1,1)$ we obtain $h_2(1)=1$. Applying Equation (\ref{eq:additivity}) to $x$ and $1$ gives 
$h_1(x)= h_1(x)\cdot h_2(1)=h_0(x\cdot 1)=h_0(x)$.
Hence $h_0=h_1=h_2$,  as claimed.
\end{proof}

\begin{theorem}\label{cor: additivemaps}
A map $f:E\rightarrow E$  is an endomorphism of $\mathcal M$ if and only if $f(x,y)= (a\cdot h(x), h(y)\cdot b)$ for all $x,y\in G$ where $h$ is an endomorphism  of  $G$  $G$ and $(a,b)\in E$.
\end{theorem}
\begin{proof}
Due to the associativity of $\cdot$,  if  $f$ has the form given above, then it is  a homomorphism of $\mathcal M$. In particular, for every $u:=(a,b)\in E$,  the transformation  $t_{u}$ defined on  $E$  by setting $t_{u}(x,y):=(a \cdot x, y\cdot b)$ is an   endomorphism of $\mathcal M$. Now, let $f$ be an endomorphism of $\mathcal M$. Let $(a,b):= f(1,1)$ and  $(a,b)^{-1}:=(a^{-1}, b^{-1})$. Since these transformations  are endomorphisms  of $\mathcal M$, the map $g:= t_{(a,b)^{-1}}\circ f$ is an endomorphism  of $\mathcal M$ fixing  $(1,1)$. According to  Lemma \ref{lem:additive}, $g$ is of the form $(h,h)$ where $h$ is a homomorphism of $G$. Since $f= t_{(a,b)}\circ g$, $f$ has the form given above.  \end{proof}

%
%
%
%
%
%
\subsection{Quasigroups,  $3$-nets and orthogonal systems}\label{quasigroup}
A \emph{quasigroup} is a magma $(G, \cdot)$ of which each translation is bijective (each element $a\in G$ yields the left translation $x \hookrightarrow a\cdot x$ and the right translation $x \hookrightarrow x\cdot a$ ). 
Let $G$ be  a quasigroup. On $E:= G\times G$, the system $\mathcal M:=(E, (\simeq_i)_{i<3})$ defined as in the case of groups is a $3$-net. The first part of Lemma \ref{lem:additive} applies and, if $G$ has an identity, the second part too. A  quasigroup $(G, \cdot)$ with  an identity  is a \emph{loop}. It is  known that  a semigroup $G$ satisfying  the identity $(x\cdot y)\cdot (z\cdot x)=x\cdot((y\cdot z)\cdot x)$ for all $x,y,z\in G$ is a loop, called a  \emph{Moufang loop}.

\begin{lemma} The $3$-net $\mathcal M$ on $G\times G$ associated with a Moufang loop $G$ with more than one element has proper automorphisms hence it is not semirigid. 
\end{lemma}

\begin{proof} Given an $a \not =1$,  let  $f:G\times G\rightarrow G\times G$ be defined by setting $f(x,y):= (a.x, y.a)$. This map is an endomorphism of $\mathcal M$. It is not constant since $G$ is a quasigroup. It is not the identity since $a\not  =1$.  Since it is not the identity, nor a constant map,  $\mathcal M$ is not semirigid. 
\end{proof}

In order to extend  Theorem \ref {cor: additivemaps} to Moufang loops it would suffice that the identity $(a\cdot x)\cdot (y\cdot b)=a\cdot((x\cdot y)\cdot b)$ holds, but this  amounts to associativity (replace $a$ by $1$), hence to the fact  that the Moufang loop is a group.  \\

Let us say that two equivalence relations $\rho$ and $\tau$ on a set $E$ are \emph{orthogonal} if each class of $\rho$ intersects each class of $\tau$ in  at most one element(that is $\rho\cap \tau =\Delta_E$)  and \emph{strongly orthogonal} if   each class of $\rho$ intersects each class of $\tau$ in exactly one  element.  According to our introduction,  a system of three equivalence relations is a \emph{$3$-net} if and only if  these relations are pairwise strongly orthogonal. 

The link  between the two notions is the following. 
\begin{proposition} \label{prop:extending}Let  $\mathcal M:=(E, (\rho_i)_{i<3})$ be a system of three equivalence relations.  These  relations are pairwise orthogonal  if and only if $\mathcal M$ is embeddable into a $3$-net. Moreover, if $E
$ is finite, then  there is a $3$-net of order at most $2\vert E\vert$  in which $\mathcal M$ is embeddable. 
\end{proposition}

\begin{proof}

We prove first that the orthogonality conditions are necessary. Next, we prove that they suffice by assuming  first that $E$ is finite, in which case we apply a result of Evans. If $E$ is infinite, we use the diagram method of Robinson, a basic technique of mathematical logic. Namely, we define a  theory $T$ in a first order language; from the finite case of our proposition, we get that  it is consistant, hence  by the Compactness theorem of first order logic it has some model;  it turns out that such a model is a  $3$-net extending $\mathcal M$. 
 
So suppose that  $\mathcal M$ is embeddable into a $3$-net $\mathcal M':=(E', (\rho'_i)_{i<3})$, let $A\subseteq E'$ be  the image of $E$ by an embedding. Since $\mathcal M'$ is a $3$-net, the $\rho'_i$'s are pairwise orthogonal,  hence their restrictions $\rho'_i\cap A\times A$ are pairwise orthogonal too. This is equivalent to the  that   the $\rho_i$'s are pairwise orthogonal. 

Conversely, suppose that these equivalence relations are  pairwise orthogonal. As in the proof of Theorem \ref{theo1}, let  $E_i:= \{\rho_i(x): x\in E\}$ be the set of equivalence classes of $\rho_i$ and let
 $\rho: E\rightarrow \prod_{i\in I} E_i$ be defined by setting $\rho(x):=(\rho_i(x))_{i\in
I}$. Let $e_i:=\vert E_i\vert$, $e:=Max\{e_i: i<3\}$. Extending each $E_i$ to a set $E'_i$ of cardinality $e$, the map $\rho$ becomes an embedding of $\mathcal M$ into $(\prod_{i<3} E'_i, (\sim_i)_{i<3} )$. Let $\rho': E\rightarrow E'_1\times E'_2$ be defined by setting $\rho'(x):=(\rho_1(x) ,\rho_2(x))$. Since $\rho_1$ and $\rho_2$ are orthogonal, the system $(E, (\rho_1, \rho_2))$ is reduced  hence,  according to Theorem \ref{theo1},  the   map $\rho'$  is an embedding from $(E, (\rho_1, \rho_2))$ into $(E'_1\times E'_2, (\sim_1, \sim_2))$. 

\noindent {\bf Case 1.} Assume that  $e$ is finite. Labelling  the elements of $E'_i$ by the integers from $1$ up to $e$,  the range $E'$ of $E$  by $\rho'$ appears as a subset of the grid $e\times e$ whose elements are labelled by non-negative integers from $1$ up to $e$ (corresponding to $\rho_0(x)$) and yields an incomplete  latin square. According to Evans \cite{evans}, this incomplete latin square extends to  a latin square of order $2e$. This latin square provides a $3$-net of order $2e$ in which $\mathcal M$ embeds. 

\noindent {\bf Case 2.} Assume  that $e$ is infinite. The definition of the theory $T$ goes as follows. The language consists of three binary predicate symbols $r_i$, $i<3$,  and infinitely many constant symbols $(c_a)_{a\in E}$.  Axioms of $T$ are such that every  model $\mathcal N:= (N, (r^{\mathcal N}_i)_{i<3}, (c^{\mathcal N}_a)_{a\in E})$  interpreting predicates and constants yields  a $3$-net $(N, (r^{\mathcal N}_i)_{i<3})$ extending   $\mathcal M$ through the embedding mapping each $a \in E$ to $(c_a)\in N$. For that,  we define  two sets of axioms; the first set about $r_i$, $i<3$,  expresses that   $(N, (r^{\mathcal N}_i)_{i<3})$ will form a $3$-net; the second set that the map $a\rightarrow  c^{\mathcal N}_a$ is an embedding from $\mathcal M$ into $(N, (r^{\mathcal N}_i)_{i<3})$, the axioms needed here are $ \neg (c_a= c_{a'}) $ for every  $a, a'$ distinct in $E$ and $ r_i( c_a,  c_{a'})$, resp. $ \neg  r_i( c_a,  c_{a'})$,  whenever $\rho_i(a, a')$, resp. $\neg \rho_i(a, a')$,  holds in $\mathcal M$.  Now, if $F$ is any  finite set of axioms, there is a finite subset $E_F$ of $E$ such that no constant symbols $c_a$ with $a\in E\setminus E_F$ appears among  the axioms belonging to   $F$.  Applying the finite case of this  proposition to $\mathcal M{\restriction E_F}$ we find  a model of $F$  in which all   constant symbols $ c_a$ with $a\in E\setminus E_F$ are interpreted by some arbitrary element $c\in E_F$. Since every finite set of axioms has a model, Compactness theorem of first order logic ensures that $T$ has a model.  \end{proof}

\subsection{Triangles and monogenic systems} \label{subsection:monogenic}
Let $\mathcal M:=(E, (\rho_i)_{i<3})$ be a system of three pairwise orthogonal equivalence relations. 
\begin{definitions}A \emph{triangle} $T$ in $E$ consists of  three elements $u_0, u_1,u_2$ (not necessarily distinct) of $E$ such that   $u_0\simeq_2 u_1$, $u_1\simeq_0 u_2$ and $u_2\simeq_1 u_0$. For $X \subseteq C \subseteq E$ we say that $X$   is $C$-\emph{closed} if  any triangle of $C$ with two elements in $X$ is  included in $C$. Denote by $\varphi_{C}(X)$ the intersection of $C$-closed subsets  of $C$ containing $X$ and  denote by  $\delta_{ C}(X)$ the set of $u\in C$ such that there is a triangle $\{a,b,u\}$ with  $\{a,b\}\subseteq X$. 
For $n \ge 0$ define $\delta^{(n)}_{ C}(X)$ recursively by $\delta^{(0)}_{\mathcal C}(X)=X$  and $\delta^{(n+1)}_{ C}(X):=\delta_{ C}(\delta^{(n)}_{ C}(X))$. \end{definitions}

Trivially, the collection of $C$-closed  subsets   is stable under intersection and  $\varphi_{ C}(X)=\bigcup \{\delta^{(n)}_{ C}(X): n\in \N\}$.

\begin{definitions}\label{def:monogenic} We say that a subset $X$ of $C$ \emph{generates} $C$ if $\varphi_{ C}(X)=C$. We say that $C$ is \emph{monogenic}
if some subset of $C$ with at most two elements generates $C$.
\end{definitions}

\begin{examples}\label{examples:semirigid}  If $G$ is a cyclic group (denoted additively)  and $E:=G\times G$,   then $E$ is monogenic (indeed, $E$ is generated by the pair $\{(0,0), (1,0)\}$ where $1$ generates $G$). Also for each $n \in \N$, the set $T_{n}=\{(i,j) \in \N \times \N : i+j \le n\}$ defined in the introductory section  is monogenic. Indeed, the set $X:=\{(0,0), (0,1)\}$ generates it. Similarly, $X:=\{(0,n),(1,n-1)\}$ generates each of $T_{n,2}$ and 
 $T'_{n,2}$ also defined in the introduction. However, the $8$-element set  $$U:=\{(0, 0),  (2, 0),(1, 1), (2, 1),(1, 2), (2, 2), (0, 3), (1, 3)\}$$  is not monogenic (see Figure \ref{2}).
\end{examples}

 \begin{lemma}[Triangle determinacy property]  Let  $i,j,k$ be the three elements of $\{0,1,2\}$ in an arbitrary order and $u,v,w\in E$ such that $w\simeq_iv, w\simeq_ju, u\simeq_kv$. Then any $w'\in E$ such that $w'\simeq_iv$ and  $w'\simeq_ju$ is equal to $w$.
 \end{lemma}
\begin{lemma}\label {lem:extension} Let $C$ be a subset of $E$ and  $f$, $g$ be two homomorphisms of $\mathcal M\restriction C$ into $\mathcal M$. If $f$ and $g$   coincide on a subset $X$ of $C$  then $f$ and $g$ coincide on $\varphi_{ C}(X)$.
\end{lemma}
\begin{proof}
It suffices to prove that if $f$ and $g$   coincide on a subset $X$ of $C$  they coincide on $\delta_{ C}(X)$. For that, let $c\in \delta_{ C}(X)$; our aim is to prove that $f(c)=g(c)$. Let $T:=\{a,b,c\}$ be  a triangle with  $\{a,b\}\subseteq X$. The  Triangle determinacy  applied  with $u:=f(a)=g(a)$,   $v:=f(b)=g(b)$, $w:=f(c)$ and $w':=g(c)$ yields $g(c)=f(c)$. 

\end{proof}
\begin{problem}
Are the maps preserving a finite monogenic system 
 either constant or  automorphisms? \end{problem}

According to Proposition \ref{prop:extending}, we may suppose that $\mathcal M$ is a $3$-net. More can be said if  the associated quasigroup is an  additive group, particularly if this group is a subgroup of  the additive group of a field.

\subsection{$3$-net of an abelian group} Here, we denote by $A$ an abelian  group, the addition being  denoted $+$ and the neutral element $0$, and we set $E:= A\times A$.

For each  non-empty subset $X$ of $E$,  let $D(X):=\{p_i(u)-p_i(v): u,v\in X, i\in \{1, 2\} \}$. For $\alpha\in X$, set $G(X):=\alpha +\gen{D(X)}\times \gen{D(X)}$ (=$\{\alpha+ (x,y): (x,y)\in \gen{D(X)}\times \gen{D(X)}$) where $\gen{D(X)}$ is the additive subgroup of $A$ generated by $D(X)$. Note that the choice of $\alpha$ is irrelevant.
\begin{lemma}\label{lem:integerenveloppe} For every  $X\subseteq C\subseteq E$,  $\varphi_{C}(X)$ is a subset of $G(X)$.
\end{lemma}
\begin{proof}
First observe that $X\subseteq G(X)$, next that for every subgroup $B$ of $A$ and any $\beta\in E$, $\beta + B\times B$ is $E$-closed hence $\delta_C (X)\subseteq G(X)$ and finally that $\gen{D(X)}=\gen{\delta(D(X))}$.
\end{proof}
The next lemma expresses the fact that monogenic systems can be viewed as  subsystems of $\Z\times \Z$ or $(\Z/n)\times (\Z/n)$ for some integer $n$.

\begin{lemma}\label{lem:zede} If $C\subseteq A\times A$  is monogenic then $\mathcal M \restriction  C$ is isomorphic to some system
 $M'\restriction C'$ where $C'\subseteq A'\times A'$ and $A'$ is a cyclic group.  Furthermore, if $C'$ contains the vertices of a triangle, then  we may assume that  the pair $\{(0,0), (c,0)\}$ (for some  $c\in A'$) generates $C'$.\end{lemma}
\begin{proof}
We may suppose that $C$ contains the vertices of a triangle, otherwise $C$ has at most two elements and the result is trivial.  Let $X:=\{a,b\}$ be  a subset of $C$ which generates $C$. Then there is some triangle $T:=\{a,b,c\} \subseteq C$. Clearly $\{b,c\}$ and $\{a,c\}$ generate $C$. Hence, with no loss of generality, we may suppose that $p_2(a)=p_2(b)$. Thus $\gen{D(X)}=\Z (p_1(b)-p_1(a))$ and $$G(X)=a+\Z (p_1(b)-p_1(a))\times \Z (p_1(b)-p_1(a)).$$  According to Lemma \ref{lem:integerenveloppe}, $C$ is included in $G(X)$. According to Theorem \ref{cor: additivemaps} we may translate $a$ and $b$ onto $(0,0)$ and $(p_1(b)-p_1(a), 0)$ respectively.
 \end{proof}

If $A$ is the additive group of a field, say $\K$,  multiplication of two elements $\lambda$ and  $\mu$ of $A$ will be denoted by  $\lambda\mu$. In this case,  a map  $h$ from $E$ into $E$ is a  \emph{homothety composed with a translation} if  for some $\alpha \in E$ and  some \emph{coefficient}\; $\lambda \in A$,  $h(u)=\lambda u +\alpha$   for all $u\in E$.  

Note that if $A'$ is a cyclic subgroup of $A$ and $h$ preserves $E':= A'\times A'$ then 
$\lambda \in A'$ and $\alpha\in E'$ (indeed, since $h$ preserves the additive group $E'$,  
$\alpha=h(0)\in E'$; it follows that  $\lambda u\in E'$ for every $u \in E'$; hence, for a generator $a'$ of $A'$,  $\lambda (a',a')= k(a',a')$ for some $k\in \Z$; this yields $\lambda \in A'$). We will say that $h$ is a  \emph{homothety composed with a translation on  $E'$}. 

\begin{proposition} \label{prop:monogenic} Suppose that $A$ is the additive group of a field. If a subset    $C\subseteq E:= A\times A$ is monogenic  then every homomorphism  of $\mathcal M\restriction C$ into $\mathcal M$ is the restriction to $C$ of the composition of  a homothety with a translation.  

\end{proposition}
\begin{proof} Let $f$ be a homomorphism  of $\mathcal M\restriction C$ into $\mathcal M$ and let $X:=\{a,b\}$ generates $C$. Set $a':=f(a)$ and $b':=f(b)$. First suppose that $a'=b'$. Let $w$ be the common value and $c_{w}$ be the constant map mapping $a$ and $b$ to $w$.  Since $f$ and $c_{w}$ coincide on $X$, Lemma \ref{lem:extension} ensures that they coincide on $C=\varphi_{\mathcal C}(X)$. Thus $f$ is constant. Now, assume that $a'\not=b'$. Let $\lambda\in A$ be such that ${b'-a'}=\lambda(b-a)$, set $\alpha:=a'-\lambda a$ and let $h$ be the selfmap of $A$ defined by setting $h(u):=\lambda u+\alpha$ for every $u\in A$. This map is an endomorphism of  $ \mathcal M$. It coincides with $f$ on $X$. According to  Lemma  \ref{lem:extension}, $f$ coincide with $h$ on $C=\varphi_{C}(X)$. \end{proof}

\begin{proposition} \label{prop:monogenic2} Suppose that $A$ is $\Z$ or $\Z/p\Z$ where $p$ is  prime.  If a subset    $C\subseteq E\times E$ is monogenic  then every endomorphism of $\mathcal M\restriction C$ is the restriction of the composition of a  homothety with a translation on $E$. 
\end{proposition}
\begin{proof} We may extend $A$ to the additive group of a field. Apply  Proposition \ref{prop:monogenic}\end{proof}

We recall that an abelian group $A$ is \emph{torsion-free} if $nx=0$ implies $n=0$ or $x=0$) ($n\in \N, x\in A$).  An abelian  group $A$ is \emph{divisible} if for every $a\in A$ and every positive integer $n$, the equation $n.x=a$ has a  solution. Torsion-free divisible groups are  vector spaces over the field $\Q$ of rational numbers. Also each torsion-free  divisible group extends to a field.

An element $\alpha\in E$ is a \emph{center of symmetry} for $C\subseteq E$ if for every $u\in C$, $2\alpha-u\in C$.

\begin{proposition} \label{prop:maintheo} Let $A$ be torsion-free and $C$ be  a monogenic subset of $E$ containing $(0,0)$. If each  line going through $(0,0)$  contains only finitely many elements of $C$ and if $C$ has no center of symmetry then $\mathcal M\restriction C$ is semirigid.
\end{proposition}
\begin{proof} Let $f$ be an endomorphism of $\mathcal M\restriction C$.   According to Lemma \ref{lem:zede} we may suppose that $A=\Z$. According to  Proposition \ref{prop:monogenic2}, the map $f$ is the composition of a  homothety with a translation, and  there are $\alpha\in E$ and $\lambda \in \Z$ such that $f(u)=\lambda u +\alpha$ for all $u \in E$. We show that $\alpha=0$ and $\lambda=1$. Suppose for a contradiction $\alpha\not=0$.  Let us iterate $f$. For each  integer $n\geq 1$, the iterated selfmap $f^n$ satisfies $f^{(n)}(u):= \lambda^n u+ (\sum_{k<n} \lambda ^k)\alpha$ for all $u \in E$.  Since $(0,0)\in C$, obviously $(\sum_{k<n} \lambda ^k) \alpha\in C$ for every $n\geq 1$. Hence, since all these elements  are multiples of $\alpha$,
they all lie on the line going through $(0,0)$ and $\alpha$. Since, by hypothese this line contains only finitely many elements of $C$,  either $\lambda=0$, in which case $f$ is constant,  or $\lambda=-1$,  in which case $f$ is a symmetry  whose center is $\frac{\alpha}{2}$, a case which is excluded by assumption.  Hence   $\alpha=0$ and thus $f(u)=\lambda u$. Since each  line going through $(0,0)$  contains only finitely many elements of $C$, proceeding in a similar way as above we get $\lambda \in \{1,-1\}$. Since $f$ is not constant and $(0,0)$ is not a center of symmetry of $C$, similarly we get $\lambda=1$.   In this case, $f$ is the identity on $C$, thus $\mathcal M$ is semirigid. \end{proof}

 Theorem \ref{maintheo} follows immediately from Proposition \ref{prop:maintheo}.

\begin{figure}[h]
\begin{center}
\includegraphics[width=5.4in]{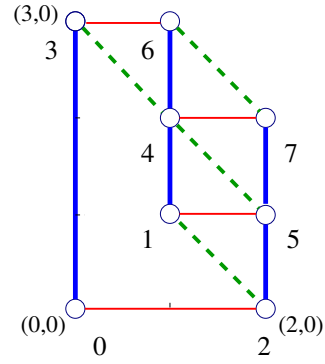}
\end{center}
\caption{An example of a not monogenic semigid system} \label{2}
\end{figure}

\begin{remark} There are non-monogenic  subsets $C$ for which $\mathcal M\restriction C$ is semirigid. Indeed,  let $A:= \Z$  and let $U$ be the ones represented Figure \ref{2}. Then $\mathcal M\restriction U$ is semirigid but $U$ is not monogenic  (the semirigidity was  checked  by computer).
\end{remark}

\subsection{Semirigid subsystems of $\R\times \R$}

Let $A:=\Z$, $E:=A\times A$ and $$B:=\{(x,y)\in \Z\times \Z:  x+y\in \{1,2\}\}\cup \{(0,0)\}.$$ 
\begin{proposition}\label {prop. semirigid}The system $\mathcal M\restriction B$ induced by the system $\mathcal M$ associated with $E$ is semirigid.
\end{proposition}
\begin{proof}  As it is easy to see, $X:=\{(0,0), (1, 0)\}$ generates $B$. Thus $B$ is monogenic. Furthermore, all other  hypothesis of  Proposition \ref{prop:maintheo} are satisfied. Hence $\mathcal M\restriction B$ is semirigid.\end{proof}

We recall that a subset $X$ of $\R$ is \emph{dense} if for every $x<y$ in $\R$ there is some $z\in X$ such that $x<z<y$. We also recall that every additive subgroup $D$ of $\R$ is either discrete, in which case $D=\Z\cdot r$ for some $r\in \R$, or dense. Let $D$ be an additive subgroup of $\R$ containing $\Z$. Set $\Delta:=\{(x,y)\in D\times D: 0\leq x, 0\leq y, x+y\leq 1\}$.

\begin{proposition}\label{prop:uncountable}
If $D$ is a dense subgroup of $\R$ including $\Q$ then the system $\mathcal M\restriction ({B\cup\Delta})$ induced by the system $\mathcal M$ associated with $D\times D$ is semirigid.
\end{proposition}

\begin{proof}
Let $g$ be an  endomorphism of $\mathcal M{\restriction {B\cup\Delta}}$. In particular, $g$ induces a homomorphism of $\mathcal M{\restriction B}$ in $\mathcal M$. Since $B$ is monogenic  then, according to Lemma \ref{prop:monogenic}, the map $g{\restriction B}$ is the restriction to $B$ of the composition of an homothety with a translation, that is   $g(u)=\lambda u+\alpha$ for all $u\in B$ and  some $\lambda \in \R$ and $\alpha \in E$. Since $g$ maps $B$ into $B\cup \Delta$, one would  easily see  that if $g$ is not constant then $\lambda=1$ and $\alpha=0$  and thus  $g$ is the identity on $B$.  Let  $\alpha:=g(0,0)$ and set $f(u):=g(u)-\alpha$ for every $u\in B\cup\Delta$. Apply Lemma \ref{lem:additive} to $E:=D\times D$, $C:=\Delta$ and to the restriction of $f$ to $C$. The maps $h_0,h_1,h_2$ 
 given by  Lemma \ref{lem:additive} coincide and their common value    $h$ satisfies $h(x+y)=h(x)+h(y)$ whenever $x,y, x+y\in [0,1]\cap D$. We claim  that $h(x)= xh(1)$ for every $x\in [0,1]\cap D$. We observe first that  $h(rx)=r h(x)$ for every rational $r$ and $x$ such that $x, rx\in [0,1]\cap D$. Next we prove that $h$ is continuous.  Our claim follows. The continuity   of $h$ will follows  from the fact that $\vert h(x)\vert \leq \vert x\vert $  for every $x\in [0,1]\cap D$ (indeed, from the additivity condition, that means  that $h$ is $1$-Lipschitz, hence continuous).  So given $x\in [0,1]\cap D$, let  $X:=[0,1]\cap \mathbb Q x=(\mathbb Q\cap [0, \frac{1}{x}])x$.
Observe that $f(X\times\{0\})=(\mathbb Q \times [0, \frac{1}{x}])h(x))\times\{y_0\}$,  thus is a set 
 $Z=X'\times\{y_0\}$ for a dense subset $X'$ of a real interval of length $\ell:=\vert\frac{h(x)}{x}\vert$.
It is easily checked that no such $Z$ can be included in $B\cup\Delta$ if $\ell>1$. Now we may conclude that either $g$ is constant or the identity. Indeed, if $g$ is constant on $B$ then, in particular  $\alpha:=g(0,0)=g(1,0)$, in which case   $h(1)=p_1(f(1,0))=p_1(g(1,0)-\alpha)=p_1(0,0)=0$. Hence, by our claim,  $h$ is $0$ on $[0,1]\cap D$. Thus  the image of $\Delta$ is $\alpha$ and $g$ is constant. If $g$ is not constant on $B$ then,  as  we have seen,  $g$ is the identity on $B$. In this case   $\alpha=(0,0)$,  hence $f=g$.  Thus $h(1)=p_1(f(1,0))= p_1(g(1,0))=p_1(1,0)=1$; then,  by our claim,  $h$ is the identity on
$[0,1]\cap D$ and  $g$ is the identity on $\Delta$. Thus,  $g$ is the identity on $B\cup\Delta$, proving that $\mathcal M{\restriction {B\cup\Delta}}$ is semirigid. \end{proof}

\subsection{Proof of Theorem \ref{thm:infinite}}
According to Z\'adori's result, there is a semirigid system of three equivalence on a set of size $\kappa$ for each finite $\kappa$ distinct from $2$ and $4$. If $\kappa=\aleph_0$, apply Proposition \ref{maintheo}. If $\aleph_0<  \kappa\leq 2^{\aleph_0}$,  observe (with the \emph{axiom of choice}) that there are dense additive subgroups of $\R$ of cardinality $\kappa$. Apply Proposition \ref{prop:uncountable}.

%
%

\noindent {\bf Acknowledgement.}

\noindent Results of this paper have been presented at  the IEEE 42nd  International Symposium on Multiple-Valued Logic, ISMVL-2012,  May 14-16, Victoria, BC, Canada,  at the special session in honor of Ivo G. Rosenberg. 
The authors are pleased to thank Pr. Vincent Gaudet, Michael Miller, Dan Simovici and the organizers of this conference. 
 They thank  Lucien Haddad for his helpful comments. 
They also thank for  Grant-in-Aid for Scientific Research (B) 21300079, 2009 (JAPAN).


\end{document}